\documentclass[preprint,12pt]{elsarticle}

\usepackage[pdftex]{hyperref}
\biboptions{sort&compress,numbers,square}

\usepackage{amsfonts}
\usepackage{graphicx}
\usepackage{epsfig}
\usepackage{geometry}
\usepackage{tikz}
\usetikzlibrary{trees,positioning,shapes}
\usepackage{subcaption}
\usepackage{amsmath}
\usepackage{amsthm}

\usepackage{xcolor, colortbl}

\usepackage{algorithm}
\usepackage{algorithmicx, algpseudocode}

\newcommand{\inred}[1]{\textcolor{red}{#1}}
\newcommand{\inblue}[1]{\textcolor{blue}{#1}}
\newcommand{\inbrown}[1]{\textcolor{brown}{#1}}

\newcommand{\inblack}[1]{\textcolor{black}{#1}}

\newcommand{\ingray}[1]{\textcolor{gray}{#1}}

\usepackage{bm}
\usepackage{mathtools}

\usepackage{rotating}

\begin{document}

\begin{frontmatter}

\title{Physics Informed Neural Networks with strong and weak residuals for advection-dominated diffusion problems.}

\author{Maciej Sikora$^{(1)}$, Patryk Krukowski$^{(1)}$, Anna Paszy\'nska$^{(2)}$, Maciej Paszy\'{n}ski$^{(1)}$}

\address{$^{(1)}$Institute of Computer Science, \\ AGH University of Science and Technology,
Krak\'{o}w, Poland \\
e-mail: maciej.paszynski@agh.edu.pl\\
$^{(2)}$ Faculty of Physics, Astronomy and Applied Computer Science, \\
Jagiellonian University, Krak\'ow, Poland \\
e-mail: anna.paszynska@uj.edu.pl}

\begin{abstract}
This paper deals with the following important research questions. Is it possible to solve challenging advection-dominated diffusion problems in one and two dimensions using Physics Informed Neural Networks (PINN) and Variational Physics Informed Neural Networks (VPINN)? How does it compare to the higher-order and continuity Finite Element Method (FEM)? How to define the loss functions for PINN and VPINN so they converge to the correct solutions? How to select points or test functions for training of PINN and VPINN?
We focus on the one-dimensional advection-dominated diffusion problem and the two-dimensional Eriksson-Johnson model problem.
We show that the standard Galerkin method for FEM cannot solve this problem.
We discuss the stabilization of the advection-dominated diffusion problem with the Petrov-Galerkin (PG) formulation and present the FEM solution obtained with the PG method.
We employ PINN and VPINN methods, defining several strong and weak loss functions.
We compare the training and solutions of PINN and VPINN methods with higher-order FEM methods.
\end{abstract}


\end{frontmatter}

\setcounter{tocdepth}{3}
\tableofcontents
\newcommand{\fintro}{f} 

\section{Introduction}

The classical way of solving PDEs numerically is based on the Finite Element Method (FEM). In FEM, we approximate the solution of the PDE by using a linear combination of the prescribed basis functions. The coefficients of the basis functions are obtained by solving a system of linear equations. The most accurate version of the FEM employs higher-order and continuity B-spline basis functions \cite{IsogeometricAnalysisProposal,IGAOvwCompImpl,IGATsplines,IzogeometrycznaMES}.

The neural networks are the universal approximators \cite{UniversalApproximators}. They can successfully replace or support the FEM computations.
Recently, there has been a growing interest in the design and training of neural networks for solving PDEs \cite{NODE,DiffEqDNN,BREVIS2021186,HP}.
The most popular method for training the DNN solutions of PDEs is Physics Informed Neural Networks (PINN) \cite{c1}. Since its introduction in 2019, there has been exponential growth in the number of papers and citations related to them (Web of Science search for "Physics Informed Neural Network"). It forms an attractive alternative for solving PDEs in comparison with traditional solvers such as the Finite Element Method (FEM) or Isogeometric Analysis (IGA).
Physics Informed Neural Network proposed in 2019 by George Karniadakis revolutionized the way in which neural networks find solutions to partial differential equations \cite{c1}
In this method, the neural network is treated as a function approximating the solution of the given partial differential equation $u (x) = PINN (x)$. The residuum of the partial differential equation and the boundary-initial conditions are assumed as the loss function. The learning process consists in sampling the loss function at different points by calculating the residuum of the PDE and the boundary conditions.
Karniadakis has also proposed Variational Physics Informed Neural Networks VPINN \cite{c3}.
VPINNs use the loss function with a weak (variational) formulation. In VPINN, we approximate the solution with a DNN (as in the PINN), but during the training process, instead of probing the loss function at points, we employ the test functions from a variational formulation to average the loss function (to average the PDE over a given domain).
In this sense, VPINN can be understood as PINN with a loss function evaluated at the quadrature points, with the distributions provided by the test functions.
Karniadakis also showed that VPINNs could be extended to $hp$-VPINNs ($hp$-Variational Physics Informed Neural Networks) \cite{c4,c5}, where by means of $hp$-adaptation ($h$-adaptation is breaking elements, and $p$-adaptation is raising the degrees of base polynomials) it is possible to solve problems with singularities.
The incorporation of the domain decomposition methods into VPINNs is also included in the RAR-PINN method \cite{c7}.
In conclusion, a family of PINN solvers based on neural networks has been developed during the last few years, ranging from PINNs, VPINNs, and hp-VPINNs to RAR-PINNs.

In this paper, we investigate the application of PINN VPINN family methods to solve the
advection-dominated diffusion, a challenging computational problem. It has several important applications, from pollution simulations \cite{Pollution} to flow and transport modeling \cite{LOS2021200}.
This problem is usually solved numerically using the finite element method.
However, for small values of ~$\varepsilon / \left\|\beta\right\|$, which means that the advection is much larger than diffusion,
it is a very difficult problem to solve.
The traditional finite element method does not work.
They encounter numerical instabilities resulting in unexpected oscillations and giving unphysical solutions.

There are several stabilization methods, from residual minimization \cite{CALO2021113214}, Streamline Upwind Petrov-Galerkin (SUPG) \cite{HUGHES198797}, and Discontinuous-Galerkin \cite{Ern} and Discontinuous Petrov-Galerkin methods \cite{https://doi.org/10.1002/num.20640}.

In this paper, we check if neural network-based methods can solve these challenging problems in a competitive way.

The new scientific contributions of our paper are the following. We show how to define the loss functions for PINN and VPINN so it can successfully solve the Eriksson-Johnson problem. We investigate strong and weak residuals as the candidates for the loss function.
We compare the PINN and VPINN methods with standard finite element method solvers.

The structure of the paper is the following. We start in Section 2 with an introduction of the finite element method for advection-dominated diffusion, and we explain why it does not work. Section 3 is devoted to the one-dimensional model problem solved with PINN and VPINN, and Section 4 is to the two-dimensional model Eriksson-Johnson problem. We conclude the paper in Section 5.

\section{Galerkin method for one-dimensional advection-dominated diffusion model problem}
We focus on the following model formulation of the advection-dominated diffusion problem in one dimension.
Find $u\in C^2(0,1)$:
\begin{eqnarray}
\underbrace{-\epsilon \frac{d^2u(x)}{dx^2}}_{\textrm{diffusion}=\epsilon}\underbrace{+1\frac{du(x)}{dx}}_{\textrm{advection"wind"}=1} = 0, x \in (0,1), \\
\inblue{-\epsilon \frac{du}{dx}(0)+u(0)=1.0,\textrm{ } u(1)=0}
\end{eqnarray}
The weak form is obtained by "averaging" (computing integrals) with distributions provided by test functions:

\begin{eqnarray}
{\inred{\int_0^1} -\epsilon \frac{d^2u(x)}{dx^2}\inbrown{v(x)dx}}+{\inred{\int_0^1}1\frac{du(x)}{dx}\inred{\inbrown{v(x)}dx}} = 0,\textrm{ }\forall \inbrown{v\in V} \\
{\inred{\int_0^1} \epsilon \frac{du(x)}{dx}\inbrown{\frac{dv(x)}{dx}}dx}+{\inred{\int_0^1}1\frac{du(x)}{dx}\inred{\inbrown{v(x)}dx}} + u(0)\inbrown{v(0)}= \inbrown{v(0)},\textrm{ }\forall \inbrown{v\in V}
\end{eqnarray}

\begin{figure}
\center{
\includegraphics[scale=0.2]{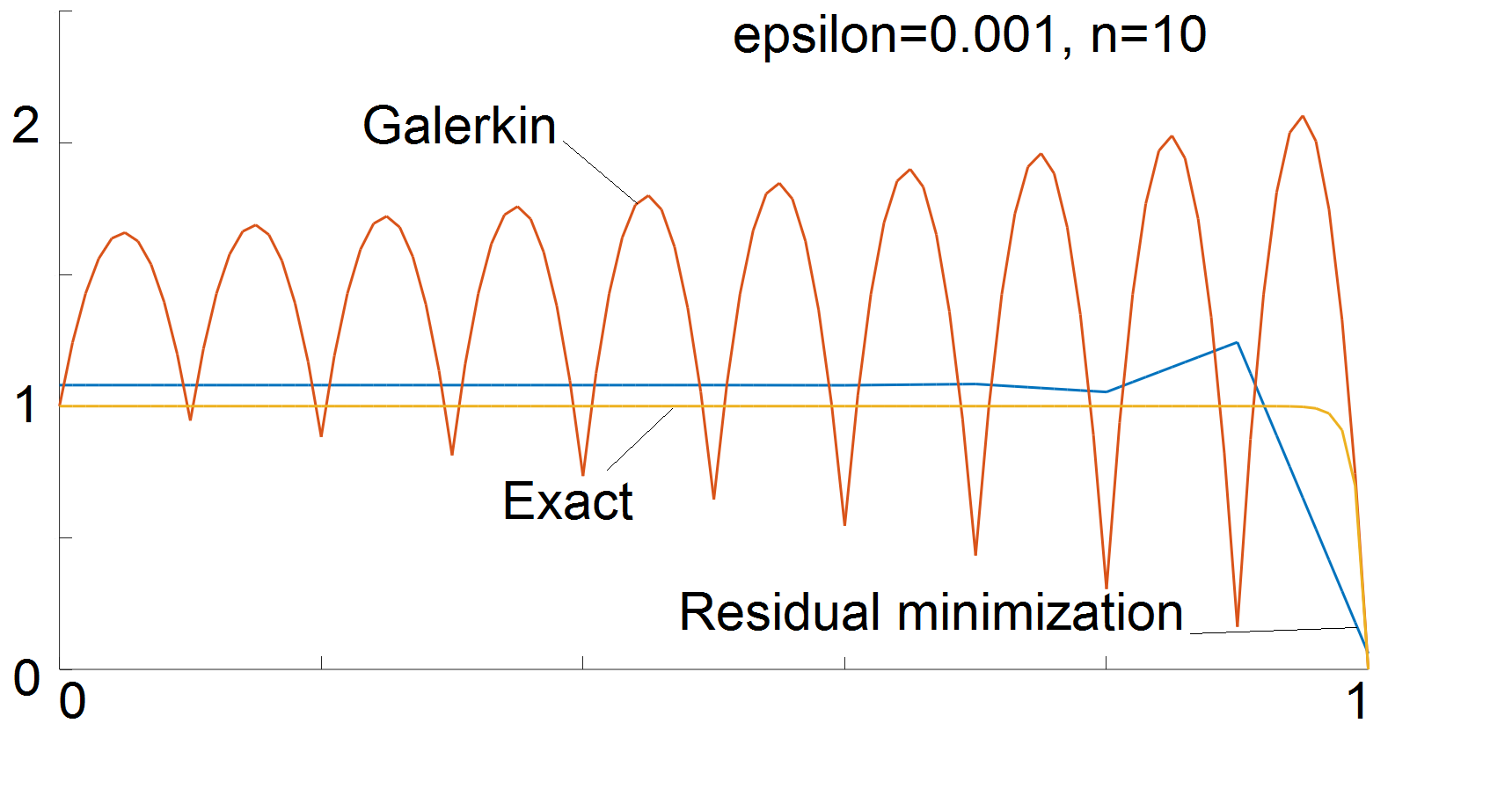} }
\caption{Comparison of the \inbrown{Galerkin method with trial=test=quadratic B-splines with $C^0$ separators}, the \inblue{Petrov-Galerkin method with linear B-splines for trial and quadratic B-splines with $C^0$ separators for test}, and the exact solution for the mesh consisting of 10 elements.}\label{fig:Fig2a}
\end{figure}

\subsection{Galerkin formulation}
\begin{figure}
\begin{centering}
\includegraphics[width=0.4\textwidth]{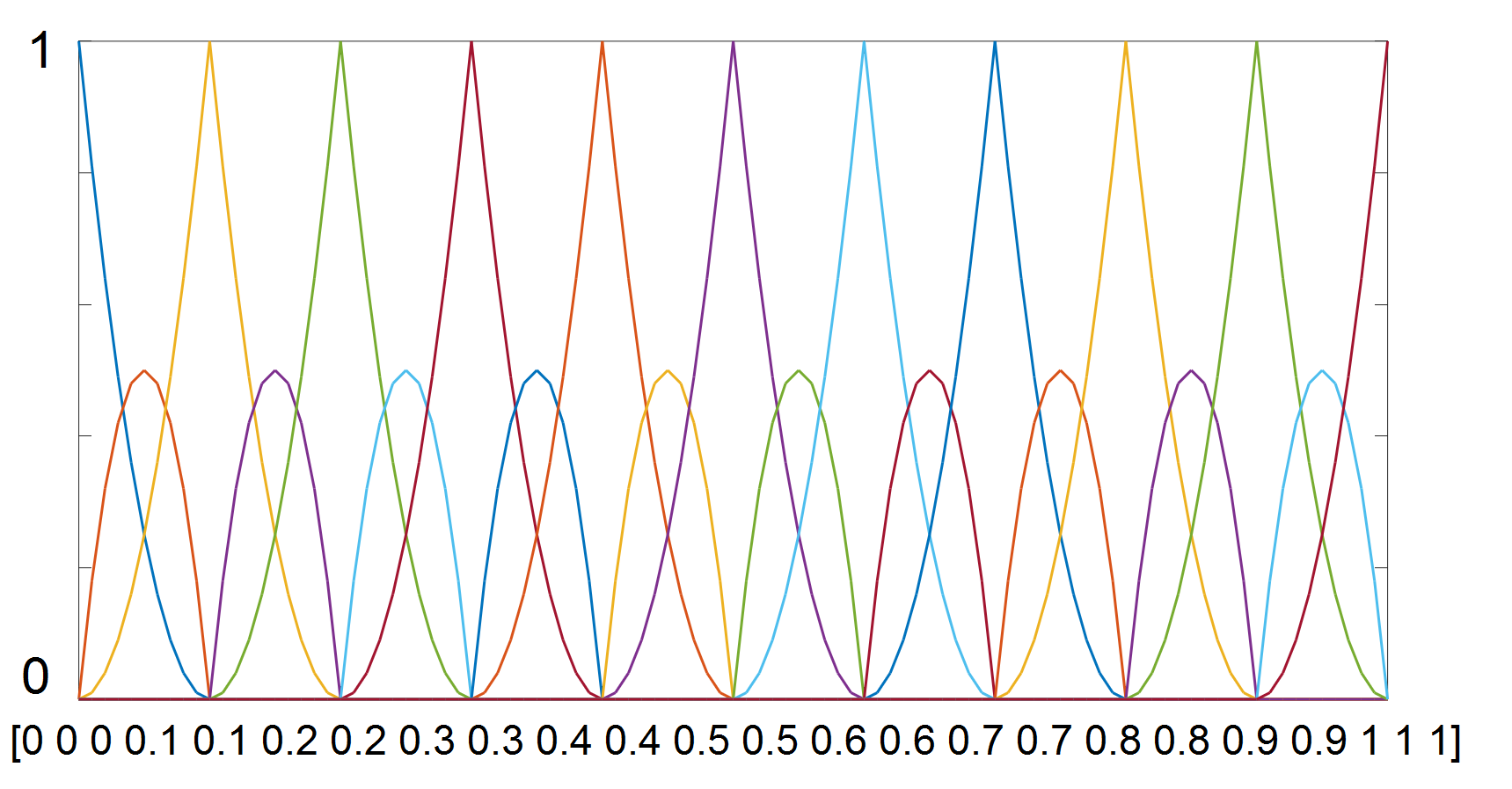}$=U_h=V_h$
\caption{Quadratic B-splines with $C^0$ separators equivalent to Lagrange basis.}
\label{fig:2}
\end{centering}
\end{figure}

The traditional finite element method is based on the Galerkin formulation, where we seek the solution as a linear combination of basis functions. In the Galerkin method, we employ the same value for trial and for testing. In our case, we can have 21 quadratic B-splines with $C^0$ separators defined by the knot vector [0 0 0 0.1 0.1 0.2 0.2 0.3 0.3 0.4 0.4 0.5 0.5 0.6 0.6 0.7 0.7 0.8 0.8 0.9 0.9 1 1 1] (see Figure \ref{fig:2}).

Find $\inblack{u_h \in U_h \subset U=V, u_h= \sum_{i=1,...,21} u_i B_i(x)}$: 
\begin{eqnarray}
{\inred{\int_0^1} \epsilon \frac{du_h(x)}{dx}\inbrown{\frac{dv_h(x)}{dx}}dx}+{\inred{\int_0^1}1\frac{du_h(x)}{dx}\inred{\inbrown{v_h(x)}dx}} + u_h(0)\inbrown{v_h(0)}= \inbrown{v_h(0)},\\ \forall \inbrown{v_h\in V_h=U_h \subset U=V}
\end{eqnarray}
where $v_h$ are the same 21 basis functions

For small values of diffusion, e.g., $\epsilon=0.001$, this problem has the following solution, presented in brown color in Figure \ref{fig:Fig2a}. As we can see from this Figure, the finite element method has problems in approximating the correct solution for the advection-dominated diffusion problem.

\subsection{Petrov-Galerkin formulation}
In the Petrov-Galerkin formulation, we employ different trial and test spaces.
\begin{figure}
\begin{centering}
\includegraphics[width=0.35\textwidth]{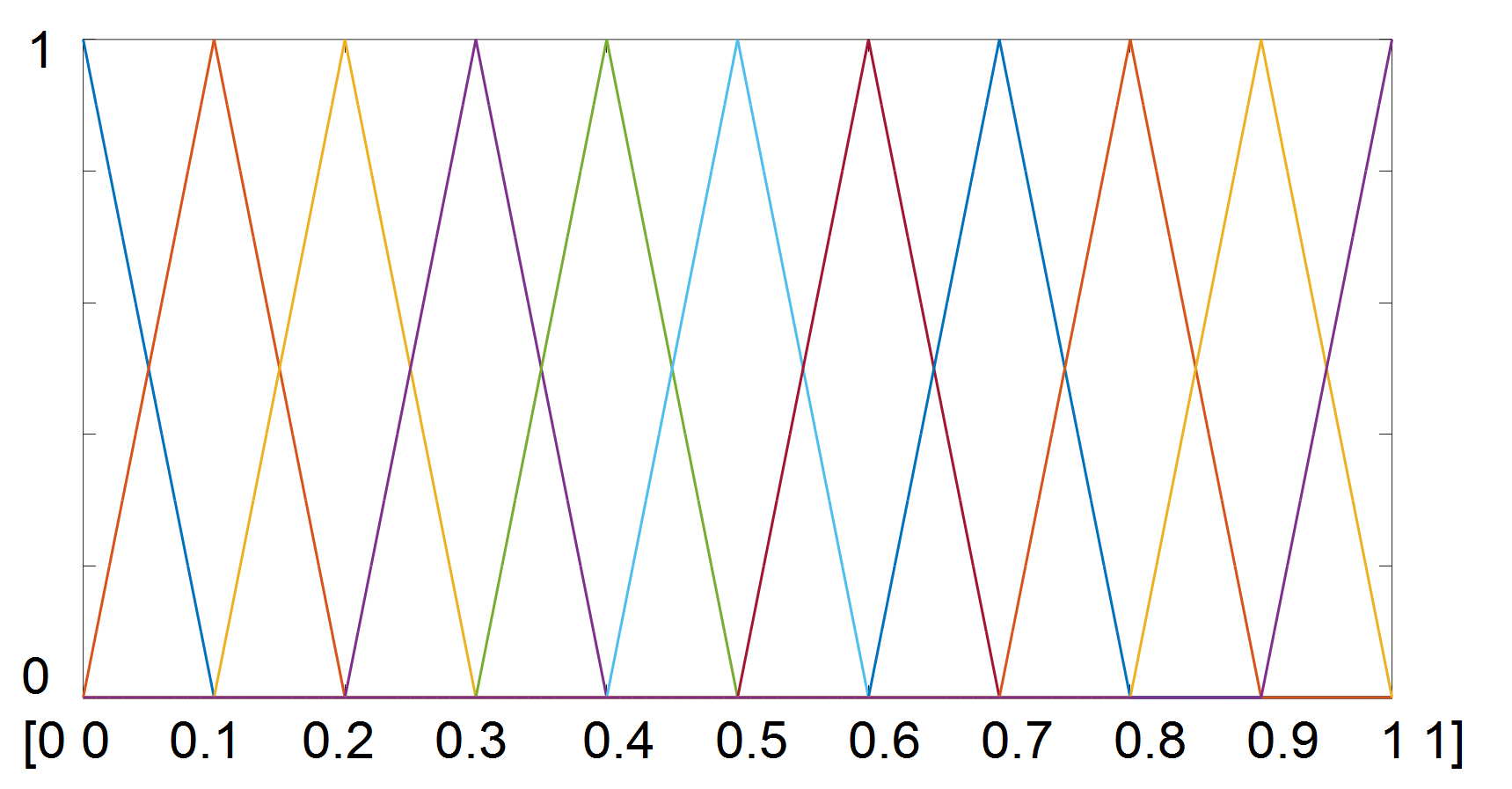}$=U_h$
\includegraphics[width=0.35\textwidth]{quadraticC0.png}$=V_h$
\caption{Linear B-splines for trial. Quadratic B-splines with $C^0$ separators equivalent to Lagrange basis for the test.}
\label{fig:3}
\end{centering}
\end{figure}

Find $\inblack{u_h \in U_h \subset U, u_h= \sum_{i=1,...,11} u_i B_i(x)}$:
\begin{eqnarray}
{\inred{\int_0^1} \epsilon \frac{du_h(x)}{dx}\inbrown{\frac{dv_h(x)}{dx}}dx}+{\inred{\int_0^1}1\frac{du_h(x)}{dx}\inred{\inbrown{v_h(x)}dx}} + u_h(0)\inbrown{v_h(0)}= \inbrown{v_h(0)},\textrm{ }\forall \inbrown{v_h\in V_h \subset V}
\end{eqnarray}
where $v_h$ are {\bf carefully selected} 11 elements of $V_h$
For example, we seek a solution as a linear combination of linear B-splines, defined with knot vector [0 0 0.1 0.2 0.3 0.4 0.5 0.6 0.7 0.8 0.9 1.0 1.0], and we test with quadratic B-splines with $C^0$ separators defined by the knot vector [0 0 0 0.1 0.1 0.2 0.2 0.3 0.3 0.4 0.4 0.5 0.5 0.6 0.6 0.7 0.7 0.8 0.8 0.9 0.9 1 1 1] (see Figure \ref{fig:3}).

\subsection{Evidence of failure of Galerkin method and superiority of Petrov-Galerkin method with optimal test (equivalent to residual minimization method) for advection-dominated diffusion problem}

\begin{figure}
\center{
\includegraphics[scale=0.2]{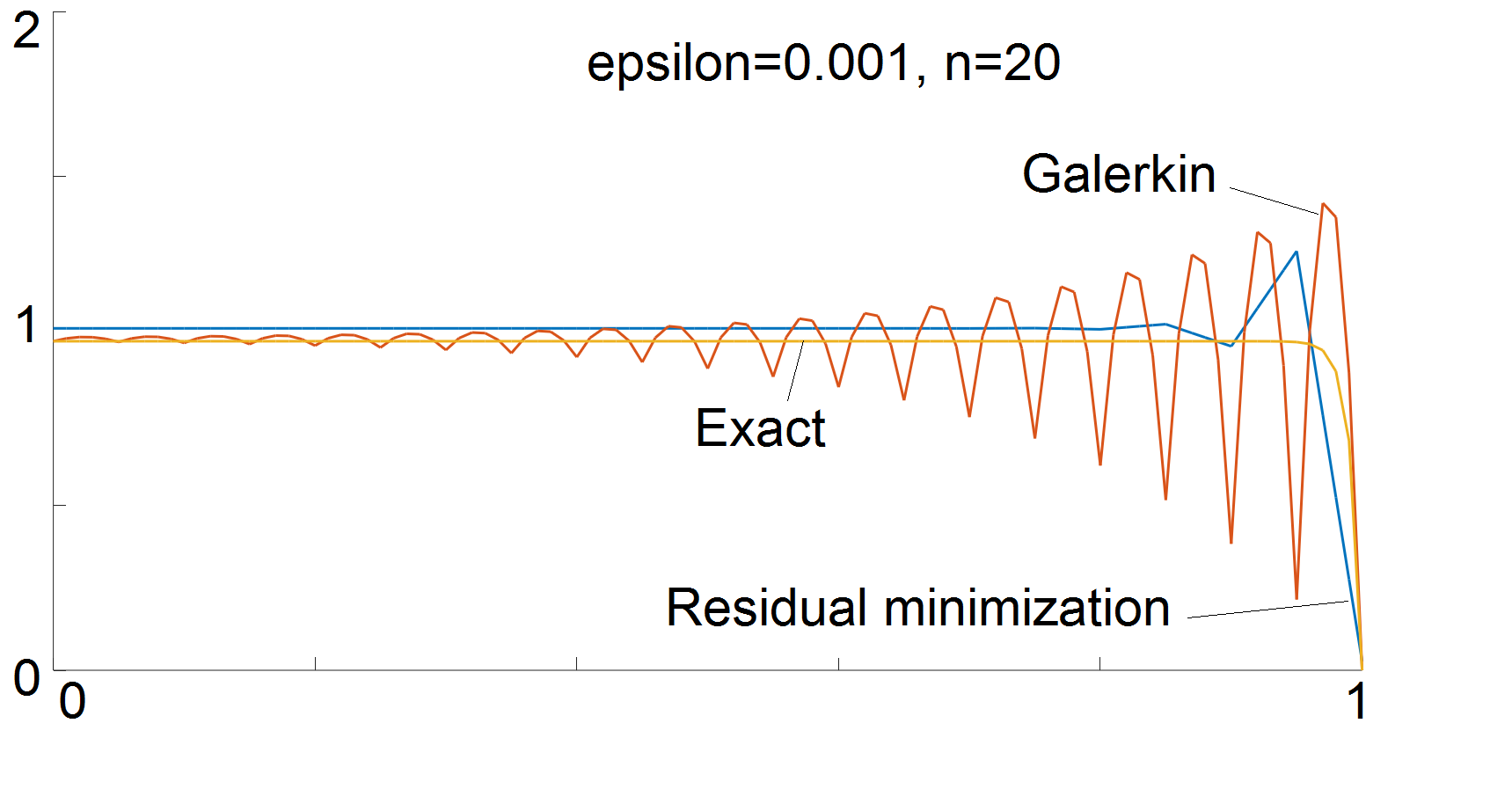} }
\caption{Comparison of the \inbrown{Galerkin method with trial=test=quadratic B-splines with $C^0$ separators}, the \inblue{Petrov-Galerkin method with linear B-splines for trial and quadratic B-splines with $C^0$ separators for test}, and the exact solution for the mesh consisting of 20 elements.}\label{fig:Fig2b}
\end{figure}

\begin{figure}
\center{
\includegraphics[scale=0.2]{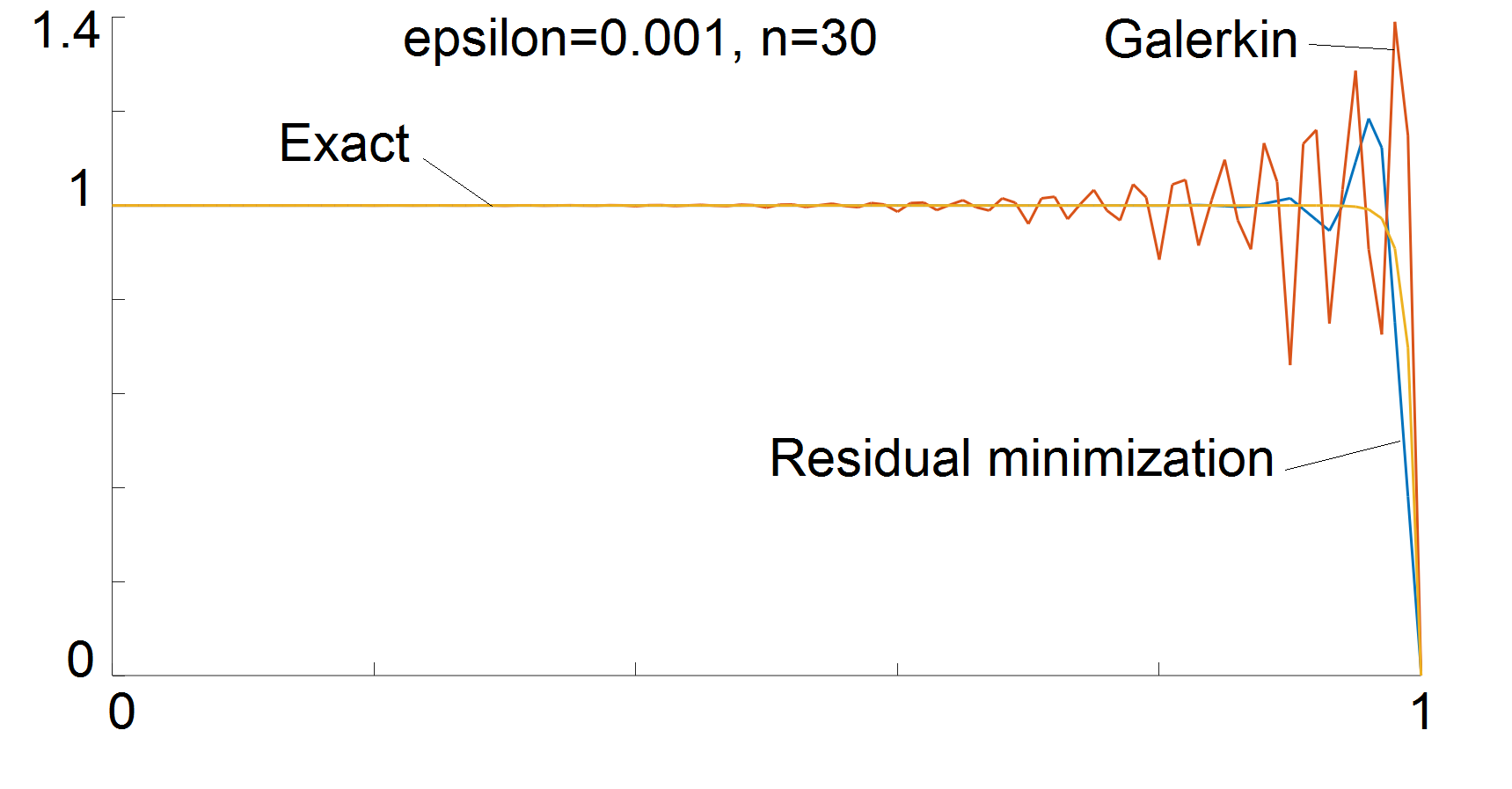} }
\caption{Comparison of the \inbrown{Galerkin method with trial=test=quadratic B-splines with $C^0$ separators}, and the \inblue{Petrov-Galerkin method with linear B-splines for trial and quadratic B-splines with $C^0$ separators for test}, and the exact solution for the mesh consisting of 30 elements.}\label{fig:Fig2c}
\end{figure}

Figures \ref{fig:Fig2a},\ref{fig:Fig2b}, and \ref{fig:Fig2c} present the comparison of the solutions for $\epsilon=0.001$ obtained with the Galerkin method using 10, 20, and 30 elements and quadratic B-splines with $C^0$ separators (equivalent to quadratic Lagrange basis) for trial and test, as well as Petrov-Galerkin method using linear B-splines for trial and quadratic B-splines with $C^0$ separators for the test.
These Figures show that the Galerkin method has strong oscillations, while the Petrov-Galerkin method stabilizes the problem (but still delivers some oscillations).

Why does the Galerkin method not work, and why does the Petrov-Galerkin method provide better (but not ideal) solutions?
We must recall the Céa Lemma, proposed in 1964 in the Ph.D. thesis of Jean Céa ''Approximation variationnelle des problèmes aux limites". The Céa lemma says that
\begin{eqnarray}
|| u - u_h || \leq \frac{M}{\alpha} dist\{\inred{U_h},u\} \notag
\end{eqnarray}
where $M$ stands for the continuity constant $b(u,v)\leq M \|u\| \|v\|$,
and $\alpha$ stands for the coercivity constant $\alpha \|u\|^2 \leq b(u,v)$.
What does it mean?
Ideally, the approximation of a solution should be as good as
the distance of the space where it lives (trial space) to the exact solution.
But due to the Céa Lemma, this is true if and only if $\frac{M}{\alpha}=1$.
We have $\alpha \|u\|^2 \leq b(u,u) \leq \|u\|^2$ which implies $\frac{M}{\alpha} \geq 1$.

The reason why $\frac{M}{\alpha}>1$ is explained by the discrete Babu\v{s}ka-Brezzi inf-sup condition \cite{BB}. The coercivity constant $\alpha$ is estimated for the abstract infinite dimensional space $V$ with the supremum over the infinite-dimensional test space. The problem is that in numerical experiments, we consider finite-dimensional spaces $u_h\in U_h, \inbrown{v_h \in V_h}$, and we compute the supremum over the finite-dimensional test space, and this supremum is smaller than infinite-dimensional one
\begin{eqnarray}
\begin{aligned}
inf_{u \in U }sup_{\inbrown{v \in V}} \frac{b(u,v)}{\|v\|\|u\|} = &\inred{\alpha} > \inred{\alpha_h} = inf_{u_h \in U_h }sup_{\inbrown{v_h \in V_h}} \frac{b(u_h,\inbrown{v_h})}{\|\inbrown{v_h}\|\|u_h\|} \notag \\ \notag &\frac{M}{\inred{\alpha_h}} > \frac{M}{\inred{\alpha}}\geq 1 \notag
\end{aligned}
\end{eqnarray}
To improve the solution, we need to test with larger test space $\inbrown{V_h}$, so it realizes supremum for $\alpha$.
This is the idea of the Petrov-Galerkin formulation. But still, in the Petrov-Galerkin formulation, the test space is still finite-dimensional, so the results may still not be ideal.

\section{Neural network approach for one-dimensional advection-dominated diffusion model problem}

\subsection{Physics Informed Neural Networks for advection-dominated diffusion}

We now introduce the PINN formulation of the one-dimensional model advection-dominated diffusion problem. We recall the strong form of the PDE: Find $u\in C^2(0,1)$:
\begin{eqnarray}
\underbrace{-\epsilon \frac{d^2u(x)}{dx^2}}_{\textrm{diffusion}=\epsilon}\underbrace{+1\frac{du(x)}{dx}}_{\textrm{advection"wind"}=1} = 0, x \in (0,1),
\inblue{-\epsilon \frac{du}{dx}(0)+u(0)=1.0,\textrm{ } u(1)=0}
\end{eqnarray}

The neural network represents the solution
\begin{equation}
u(x)=\inred{NN(x)}=A_n \sigma\left(A_{n-1}\sigma(...\sigma(A_1x+B_1)...+B_{n-1}\right)+B_n
\end{equation}

where $\sigma$ is the activation faction, e.g., sigmoid $\sigma(x)=\frac{1}{1+exp(-x)}$.
We define the following loss functions related to the point-wise residual of the PDE and with the boundary conditions:
\begin{eqnarray}
LOSS_{PDE}(x) = (-\epsilon \frac{d^2\inred{NN(x)}}{dx^2}+\frac{d\inred{NN(x)}}{dx})^2, \\
LOSS_{BC0} = ( -\epsilon \frac{d\inred{NN(0)}}{dx}+\inred{NN(0)}-1.0)^2, \qquad
LOSS_{BC1} = (\inred{NN(1)}s)^2, \\
LOSS(x) = LOSS_{PDE}(x)+LOSS_{BC0}+LOSS_{BC1}.
\end{eqnarray}

The argument to the loss functions $LOSS$ is the point $x$ selected during the training process.

\subsection{Variational Physics Informed Neural Networks for advection-dominated diffusion problem}

We also introduce the VPINN formulation of the one-dimensional model advection-dominated diffusion problem. We recall the weak form of the PDE: Find $u\in H^1(0,1)$:
\begin{eqnarray}
\inred{ \int_0^1 \epsilon \frac{d\inblack{u(x)}} {dx} \frac{d\inbrown{v(x)}}{dx}dx+\int_0^1\frac{d\inblack{u(x)}}{dx}\inbrown{v(x)}dx}
+u(0)\inbrown{v(0)}
=\inbrown{v(0)}
\forall \inbrown{v\in V} \notag
\end{eqnarray}
Now, \inred{neural network IS the solution }
\begin{equation}
u(x)=\inred{NN(x)}=A_n \sigma\left(A_{n-1}\sigma(...\sigma(A_1x+B_1)...+B_{n-1}\right)+B_n
\end{equation}
For the VPINN method, we can define the following two alternative loss functions.
The first one is related to the strong residual of the PDE, multiplied by the test functions, without integration by parts, and with the boundary conditions:
\begin{eqnarray}
b_{strong}(v) = \int \left(-\epsilon\frac{d^2\inred{NN(x)}}{dx^2}v(x)+\frac{d\inred{NN(x)}}{dx}v(x) \right)dx; \quad l_{strong}(v)=0 \\
LOSS_{strong}(v) = \left(b_{strong}(v)-l_{strong}(v)\right)^2, \\
LOSS_{BC0} = (-\epsilon \frac{\inred{dNN(0)}}{dx}+\inred{NN(0)}-1.0)^2, \quad
LOSS_{BC1} = (\inred{NN(1)})^2, \\
LOSS_1(v) = LOSS_{strong}(v)+LOSS_{BC0}+LOSS_{BC1}.
\end{eqnarray}
The second one is with the weak residual of the PDE, after integration by parts, and with the boundary conditions:
\begin{eqnarray}
b_{weak}(v) = \int \left(\epsilon\frac{d\inred{NN(x)}}{dx}\frac{dv(x)}{dx}+\frac{d\inred{NN(x)}}{dx}v(x) \right)dx+\inred{NN(0)}\inbrown{v(0)}, \\ l_{weak}(v)=\inbrown{v(0)}, \\
LOSS_{weak}(v) = \left(b_{weak}(v)-l_{weak}(v)\right)^2, \\
LOSS_{BC0} = (-\epsilon \frac{\inred{dNN(0)}}{dx}+\inred{NN(0)}-1.0)^2, \quad
LOSS_{BC1} = (\inred{NN(1)})^2, \\
LOSS_2(v) = LOSS_{weak}(v)+LOSS_{BC0}+LOSS_{BC1}.
\end{eqnarray}

The argument to the loss functions $LOSS_{weak}$ and $LOSS_{strong}$ is now the test function $v$. In this paper, we select the cubic B-splines during the training process.

\subsection{Numerical results}

In this section, we summarize our numerical experiments performed for the one-dimensional advection-dominated diffusion problems with PINN and VPINN methods.

\subsubsection{PINN and VPINN on uniform mesh}

\begin{figure}[h]
\includegraphics[scale=0.45]{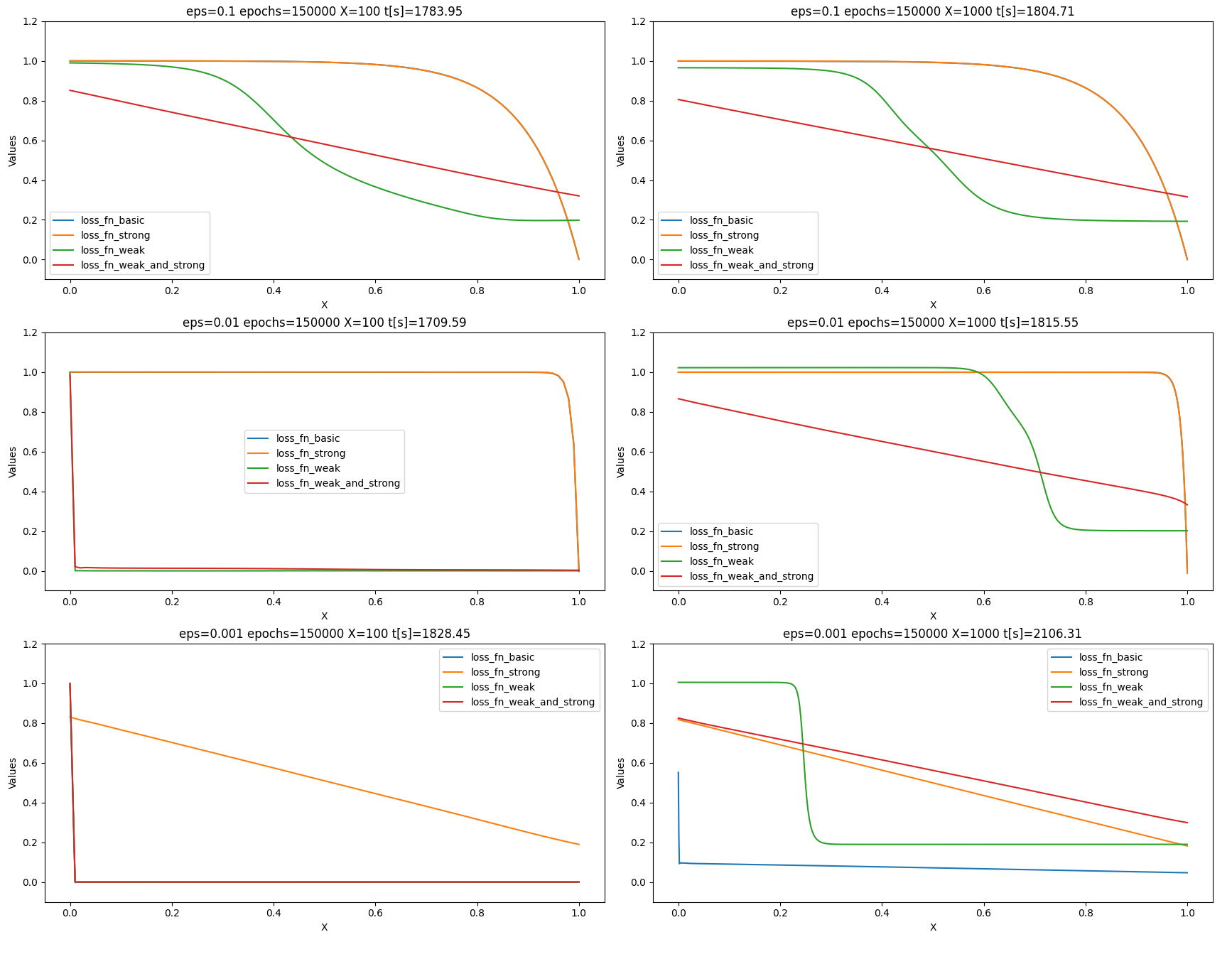}
\caption{Solution profiles of PINN and VPINN on uniform mesh}
\label{fig:solution_uniform_1D}
\end{figure}

\begin{figure}[h]
\includegraphics[scale=0.45]{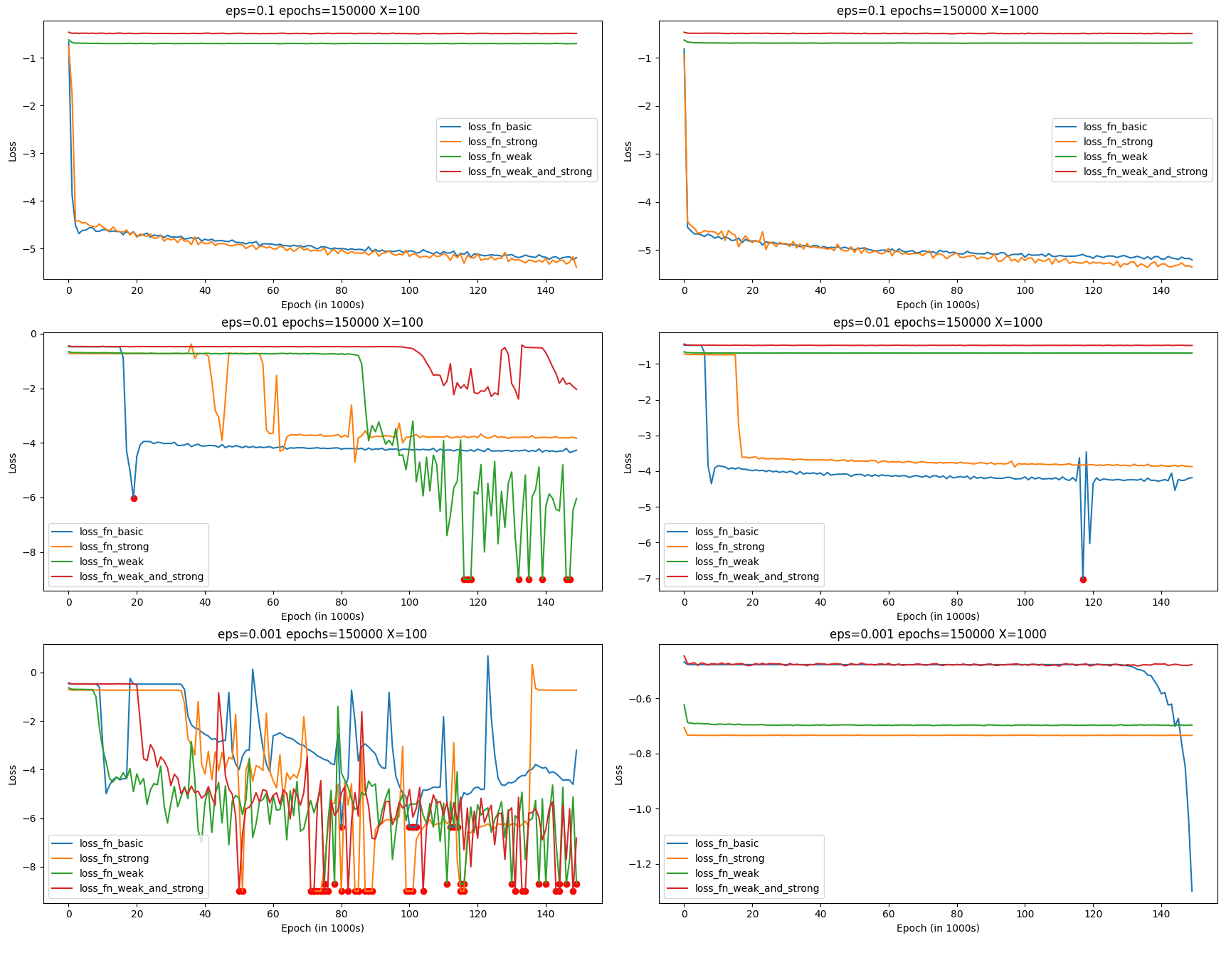}
\caption{Convergence of PINN and VPINN on uniform mesh}
\label{fig:convergence_uniform_1D}
\end{figure}

The first numerical experiment concerns the PINN method with points $x$ selected for the training with a uniform distribution and with cubic B-splines selected for VPINN using uniformly distributed intervals spanning over these points.
We employ 4 layers with 20 neurons per layer, hyperbolic tangent as the activation function, and the learning rate is defined as $\eta = 0.00125$.
We use Adam optimizer \cite{kingma2017adam}.
The summary of the solution plots is presented in Figure \ref{fig:solution_uniform_1D}. The convergence of the method is summarized in Figure \ref{fig:convergence_uniform_1D}. There are the following plots on these Figures:
\begin{itemize}
\item \emph{basic loss} corresponding to PINN method
\item \emph{strong loss} corresponding to VPINN method using strong formulation multiplied by the test functions
\item \emph{weak loss} corresponding to VPINN method using weak formulation multiplied by the test functions and integrated by parts
\item \emph{strong and weak loss} corresponding to VPINN method using both strong and weak formulations
\end{itemize}
The rows correspond to different values of $\varepsilon=0.1,0.01,0.001$, and columns correspond to a different number of points (or intervals), $X=100,1000$.
We can read from these Figures that if PINN and VPINN were using a uniform distribution of points, then it is possible to solve the problem for $\epsilon=0.01$ using 100 points.
We can conclude that for 100 points, the VPINN method with strong form loss function and the PINN method provide correct solutions for $\varepsilon=0.01$. There are 150,000 epochs, and the total training time is 1709 seconds. It is also not possible to solve the problem on uniform mesh for $\epsilon=0.001$.

\clearpage
\subsubsection{PINN and VPINN on adaptive mesh}

\begin{figure}[h]
\includegraphics[scale=0.45]{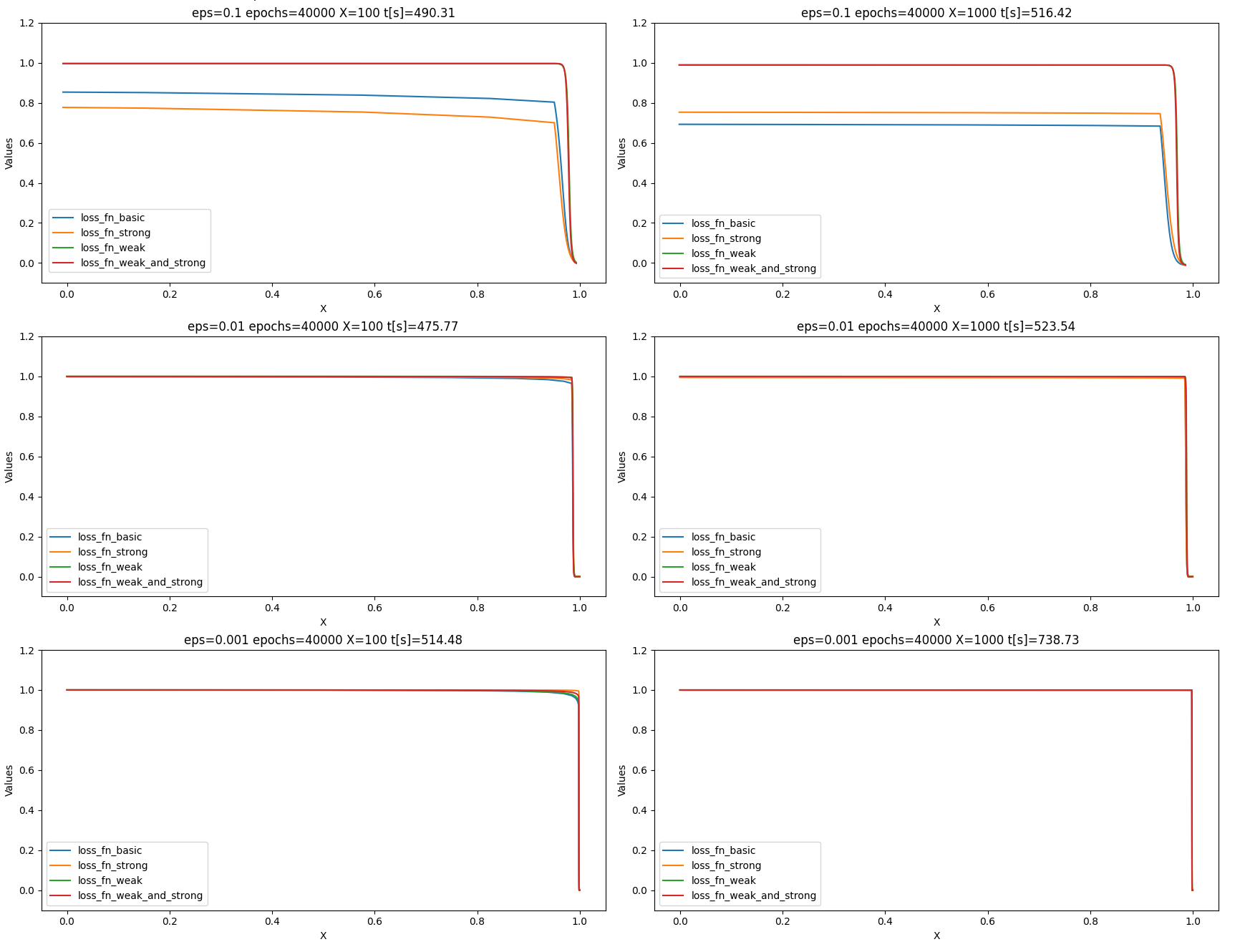}
\caption{Solution profiles of PINN and VPINN on adaptive mesh}
\label{fig:solution_adaptive_1D}
\end{figure}

\begin{figure}[h]
\includegraphics[scale=0.45]{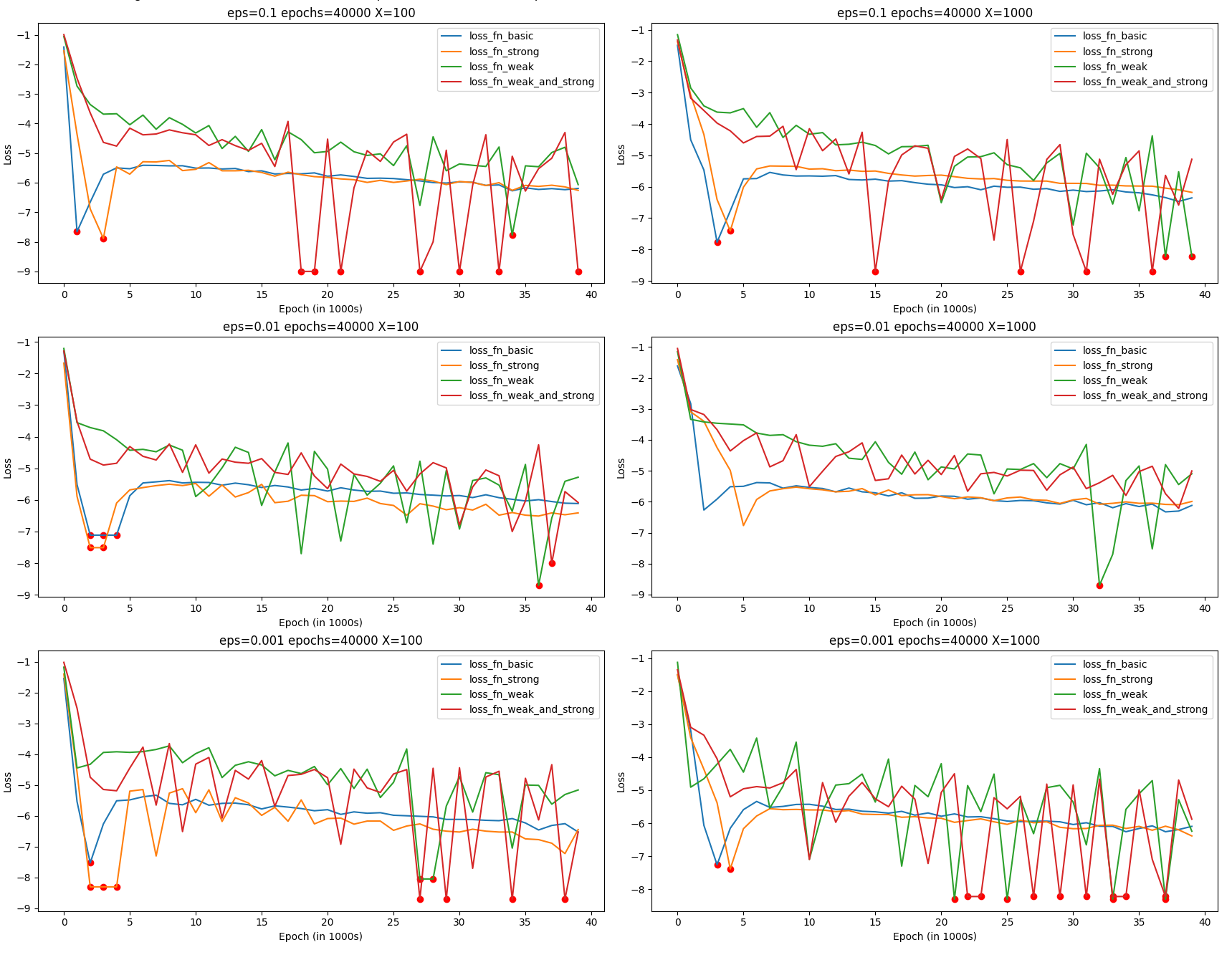}
\caption{Convergence of PINN and VPINN on adaptive mesh}
\label{fig:convergence_adaptive_1D}
\end{figure}

The second numerical experiment concerns the PINN method with points $x$ selected for the training based on the adapted mesh and with cubic B-splines selected for the VPINN span on the adapted mesh.
The mesh has been defined as $x_0=0, x_1=0.5, x_i=x_{i-1}+\frac{x_{i-1}+x_{i-2}}{2}$, up to the point where $1-x_i<\epsilon$ and then we put equally distributed remainning points between $1-\epsilon$ and 1.
The summary of the solution plots is presented in Figure \ref{fig:solution_adaptive_1D}. The convergence of the method is summarized in Figure \ref{fig:convergence_adaptive_1D}. There are the following plots on these Figures:
\begin{itemize}
\item \emph{basic loss} corresponding to PINN method
\item \emph{strong loss} corresponding to VPINN method using strong formulation multiplied by the test functions
\item \emph{weak loss} corresponding to VPINN method using weak formulation multiplied by the test functions and integrated by parts
\item \emph{strong and weak loss} corresponding to VPINN method using both strong and weak formulations
\end{itemize}
The rows correspond to different values of $\varepsilon=0.1,0.01,0.001$, and columns correspond to different numbers of points (or intervales), $X=100,1000$.
We can read from these Figures that if using adapted mesh, it is possible to solve the problems with arbitrary $\epsilon=0.1,0.01$ and $0.001$ using 100 or more points.
We can conclude that for 100 adaptively distributed points, all PINN and VPINN methods provide the correct solution for $\varepsilon=0.001$. There are 40,000 epochs, and the average training time is 514 seconds.

\clearpage

\section{Galerkin method for two-dimensional advection-dominated diffusion model problem}

This section compares the application of Petrov-Galerkin and neural-network-based methods for the solution of two-dimensional advection-dominated diffusion problems. We focus our attention on the Eriksson-Johnson model problem introduced in \cite{a3855ab1-24a7-34e2-9b85-f7a1af6ed2ff}.

\subsection{Eriksson-Johnson model problem}
Given $\Omega=(0,1)^2$, $\beta=(\beta_x,\beta_y)=(1,0)$, we seek the solution of the advection-diffusion problem\begin{equation}
\beta_x\frac{\partial u}{\partial x}+\beta_y\frac{\partial u}{\partial x}-\epsilon \left(\frac{\partial^2 u}{\partial x^2}+
\frac{\partial^2 u}{\partial y^2}\right)=0 \nonumber
\label{eq:Erikkson}
\end{equation}
with Dirichlet boundary conditions
$$
u(x,y)=0 \textrm{ for }x\in(0,1),y\in\{0,1\} \qquad
u(x.y)=g(x,y)=sin(\Pi y) \textrm{ for }x=0
$$
The problem is driven by the inflow Dirichlet boundary condition.
It develops a boundary layer of width $\epsilon$ at the outflow $x = 1$.

\begin{figure}[h]
\begin{centering}
\includegraphics[scale=0.2]{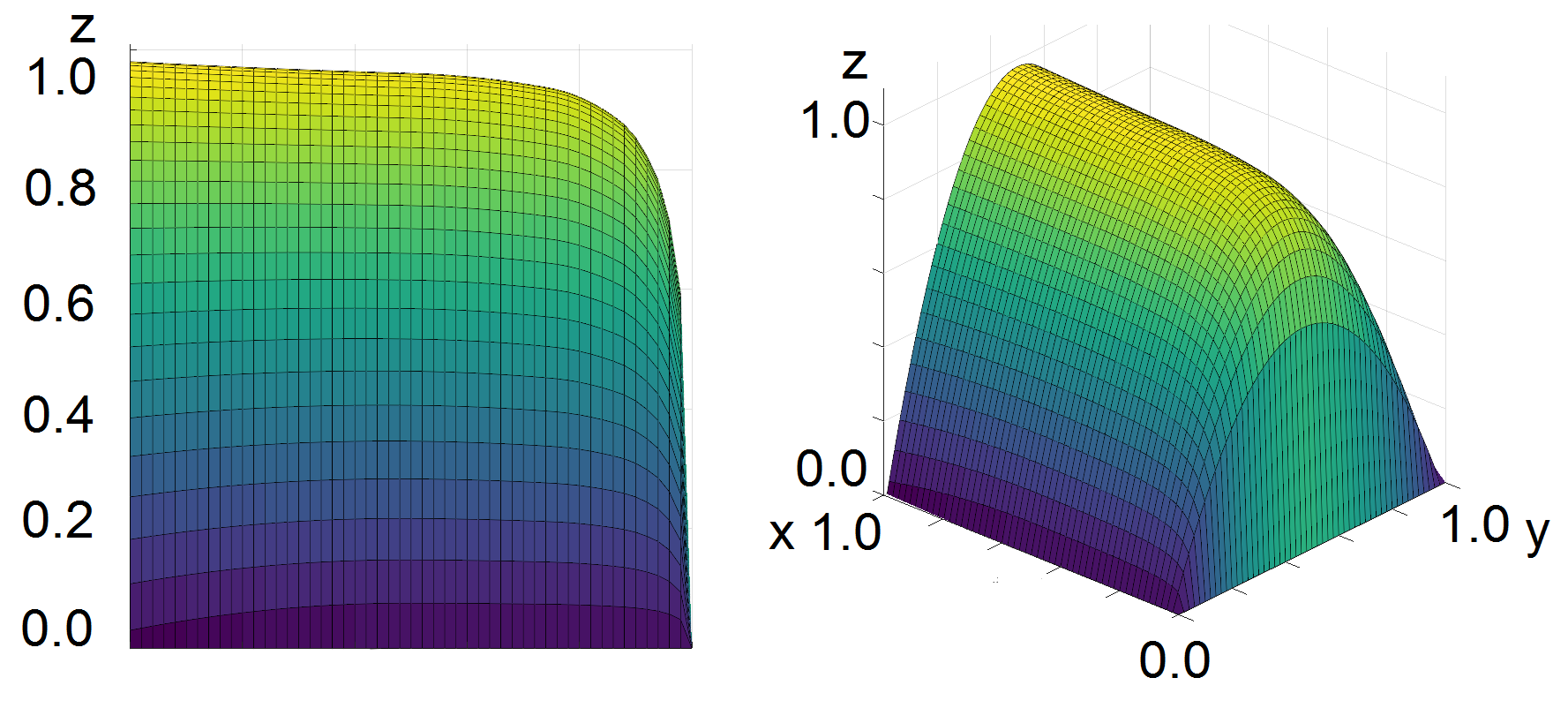}
\caption{Solution to the Eriksson-Johnson problem for $\epsilon=0.1$.}
\end{centering}
\label{fig:Erikkson_problem}
\end{figure}

\subsection{Eriksson-Johnson problem weak form for residual minimization method}

Following \cite{CALO2021113214}, we develop the weak formulation for the Eriksson-Johnson problem, with weak enforcement of the Dirichlet boundary condition.
\begin{eqnarray}
\begin{aligned}
b(u,v) &=
\beta_x\left(\frac{\partial u}{\partial x},v\right)_{\Omega}+\beta_y\left(\frac{\partial u}{\partial y},v\right)_{\Omega}+\epsilon \left( \frac{\partial u}{\partial x}, \frac{\partial v}{\partial x}\right)_{\Omega} +\epsilon \left( \frac{\partial u }{\partial y}, \frac{\partial v}{\partial y}\right)_{\Omega} \\
&\inbrown{-\left(\epsilon\frac{\partial u}{\partial x}n_x,v\right)_{\Gamma}
-\left(\epsilon\frac{\partial u}{\partial y}n_y,v\right)_{\Gamma}} \nonumber\\
&\inred{-\left(u,\epsilon \nabla v \cdot n\right)_{\Gamma} +\left(u,\beta\cdot n v\right)_{\Gamma}}\ingray{-\sum_K\left(u,3p^2 \epsilon/{h_K} v\right){_{\Gamma|_K}} } \nonumber
\label{eq:ModelProblem_b}
\end{aligned}
\end{eqnarray}
where $n=(n_x,n_y)$ is the versor normal to $\Gamma$, and $h_K$ is element diameter
\begin{equation}
l(v)=\inred{-\left(g,\epsilon\nabla v \cdot n\right)_{\Gamma}
+\left(g,\beta\cdot n v\right)_{\Gamma}}\ingray{-\left(g,3p^2 \epsilon/h v\right)_{\Gamma}}
\label{eq:ModelProblem_l}
\end{equation}
where the \ingray{gray} and \inred{red} represents the penalty terms, while the \inbrown{brown} terms result from the integration by parts.


In our Erikkson-Johnson problem, we seek a solution in space
$U = V = H^1\left(\Omega\right)$.
The inner product in $V$ is defined as
\begin{equation}
\left(u,v\right)_V=\left(u,v\right)_{L_2}+\left(\frac{\partial u}{\partial x},\frac{\partial v}{\partial x}\right)_{L_2}
+\left(\frac{\partial u}{\partial y},\frac{\partial v}{\partial y}\right)_{L_2}
\label{eq:Erikkson_inner}
\end{equation}

Keeping in mind our definitions of the bilinear form (\ref{eq:ModelProblem_b}), right-hand-side (\ref{eq:ModelProblem_l}), and the inner product (\ref{eq:Erikkson_inner}), we solve the Eriksson-Johnson problem using the residual minimization method:

Find $(r_m,u_h)_{V_m\times U_h}$ such as
\begin{eqnarray}
(r_m,v_m)_{V_m} - \left( \frac {\partial u_h}{\partial x}, v_m\right) - \epsilon\left(\frac{\partial u_h}{\partial x}, \frac{\partial v_m}{\partial x}
+\frac{\partial u_h}{\partial y}, \frac{\partial v_m}{\partial y} \right) - \notag \\
\inbrown{-\left(\epsilon\frac{\partial u_h}{\partial x}n_x,v_m\right)_{\Gamma}
-\left(\epsilon\frac{\partial u_h}{\partial y}n_y,v_m\right)_{\Gamma}} - \notag \\
\inred{-\left(u_h,\epsilon \nabla v_m \cdot n\right)_{\Gamma} +\left(u_h,\beta\cdot n v_m\right)_{\Gamma}}\ingray{-\sum_K\left(u_h,3p^2 \epsilon/{h_K} v_m\right){_{\Gamma|_K}} }
= \notag \\ \inred{-\left(g,\epsilon\nabla v_m \cdot n\right)_{\Gamma}
+\left(g,\beta\cdot n v_m\right)_{\Gamma}}\ingray{-\left(g,3p^2 \epsilon/h v_m\right)_{\Gamma}}
\quad \forall v_m \in V_h
\nonumber\\
\left( \frac {\partial w_h}{\partial x}, r_m\right) + \epsilon\left( \frac{\partial w_h}{\partial x}, \frac{\partial r_m}{\partial x}
+\frac{\partial w_h}{\partial y}, \frac{\partial r_m}{\partial y}\right)
\inbrown{-\left(\epsilon\frac{\partial w_h}{\partial x}n_x,r_m\right)_{\Gamma}
-\left(\epsilon\frac{\partial w_h}{\partial y}n_y,r_m\right)_{\Gamma}} - \notag \\
\inred{-\left(w_h,\epsilon \nabla r_m \cdot n\right)_{\Gamma} +\left(w_h,\beta\cdot n r_m\right)_{\Gamma}}\ingray{-\sum_K\left(w_h,3p^2 \epsilon/{h_K} r_m\right){_{\Gamma|_K}} }
= 0 \quad \forall w_h \in U_h
\label{eq:resmin}
\end{eqnarray}
where $(r_m,v_m)_{V_m}=(r_m,v_m)+(\frac{\partial r_m}{\partial x},\frac{\partial v_m}{\partial x})+(\frac{\partial r_m}{\partial y},\frac{\partial v_m}{\partial y})$ is the $H^1$ norm induced inner product.

\subsection{Eriksson-Johnson problem weak form for Streamline Upwind Petrov-Galerkin method}

The alternative stabilization technique is the SUPG method \cite{HUGHES198797}.
In this method, we modify the weak form in the following way
\begin{equation}
b(u_h,v_h) + \sum_K \inred{(R(u_h),\tau \beta\cdot \nabla v_h)_K}=l(v_h)
{+\sum_K \inblue{(f,\tau \beta\cdot \nabla v_h)_K}}
\quad \forall v\in V
\label{eq:Erikkson_weak_SUPG}
\end{equation}
where $R(u_h)=\beta \cdot \nabla u_h +\epsilon \Delta u_h$, and $\tau^{-1}= \beta \cdot \left(\frac{1}{h^x_K},\frac{1}{h^y_K} \right) + 3p^2\epsilon \frac{1}{{h^x_K}^2+h{^y_K}^2}$, where $\epsilon$ stands for the diffusion term, and $\beta = (\beta_x,\beta_y)$ for the convection term, and $h^x_K$ and $h^y_K$ are horizontal and vertical dimensions of an element $K$. Thus, we have
\begin{equation}
b_{\inred{SUPG}}(u_h,v_h)=l_{\inblue{SUPG}}(v_h) \quad \forall v_h\in V_h
\label{eq:SUPG_weak}
\end{equation}
\begin{eqnarray}
\begin{aligned}
b_{\inred{SUPG}}(u_h,v_h)=
& \left(\frac{\partial u_h}{\partial x},v\right)+\epsilon \left( \frac{\partial u_h}{\partial x}, \frac{\partial v_h}{\partial x}\right) +\epsilon \left( \frac{\partial u_h }{\partial y}, \frac{\partial v_h}{\partial y}\right) \\
&-\left(\epsilon\frac{\partial u_h}{\partial x}n_x,v_h\right)_{\Gamma}
-\left(\epsilon\frac{\partial u_h}{\partial y}n_y,v_h\right)_{\Gamma}\nonumber\\
&-\left(u_h,\epsilon \nabla v_h \cdot n\right)_{\Gamma} -\left(u_h,\beta\cdot n v_h\right)_{\Gamma}-\left(u_h,3p^2 \epsilon/h v_h\right)_{\Gamma} \nonumber
\\
&+\inred{\left(\frac{\partial u_h}{\partial x}+\epsilon \Delta u_h,
\left(\frac{1}{h_x} + 3\epsilon \frac{p^2}{{h^x_K}^2+{h^y_K}^2}\right)^{-1} \frac{\partial v_h}{\partial x}
\right)}
\label{eq:Erikkson_b_SUPG}
\end{aligned}
\end{eqnarray}
\begin{eqnarray}
\begin{aligned}
l_{\inblue{SUPG}}(v_h)&=(f,v_h)-\left(g,\epsilon\nabla v_h \cdot n\right)_{\Gamma} \nonumber \\
&+\left(g,\beta\cdot n v_h\right)_{\Gamma^{-}}-\left(g,3p^2 \epsilon/h_K v_h\right)_{\Gamma}
{+\inblue{\left(f,\left(\frac{1}{h_x} + 3\epsilon \frac{p^2}{{h^x_K}^2+{h^y_K}^2}\right)^{-1}\frac{\partial v_h}{\partial x}\right)}}
{.} \nonumber
\label{eq:ModelProblem_l}
\end{aligned}
\end{eqnarray}

\subsection{Eriksson-Johnson problem on adapted mesh}

\begin{figure}
\centering
\includegraphics[scale=0.24]{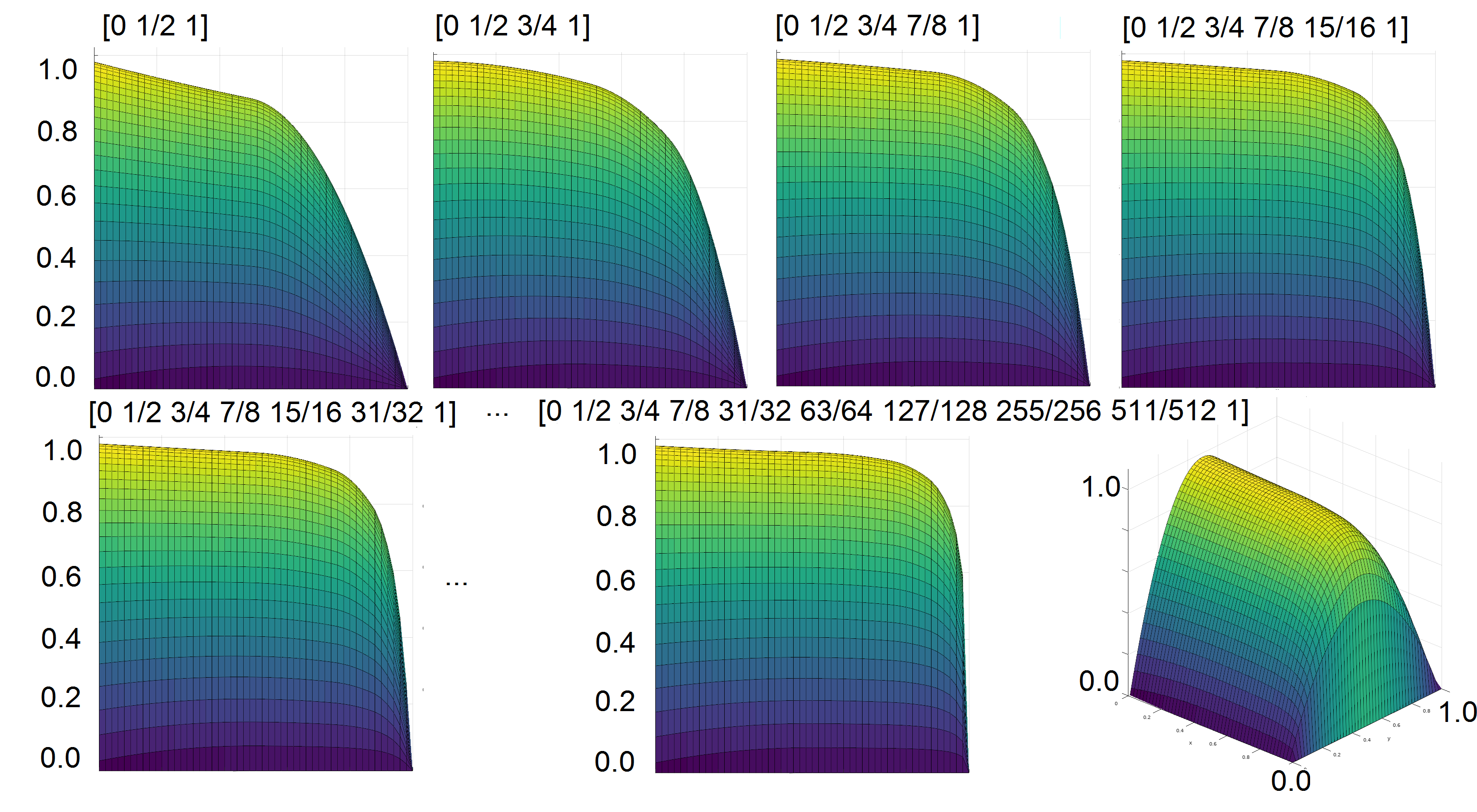}
\caption{Erikkson-Johnsson, SUPG method, $\varepsilon=10^{-3}$, using quadratic B-splines $C^1$ for trial and testing}
\label{fig:Problem3l}
\end{figure}

\begin{figure}
\centering
\includegraphics[scale=0.24]{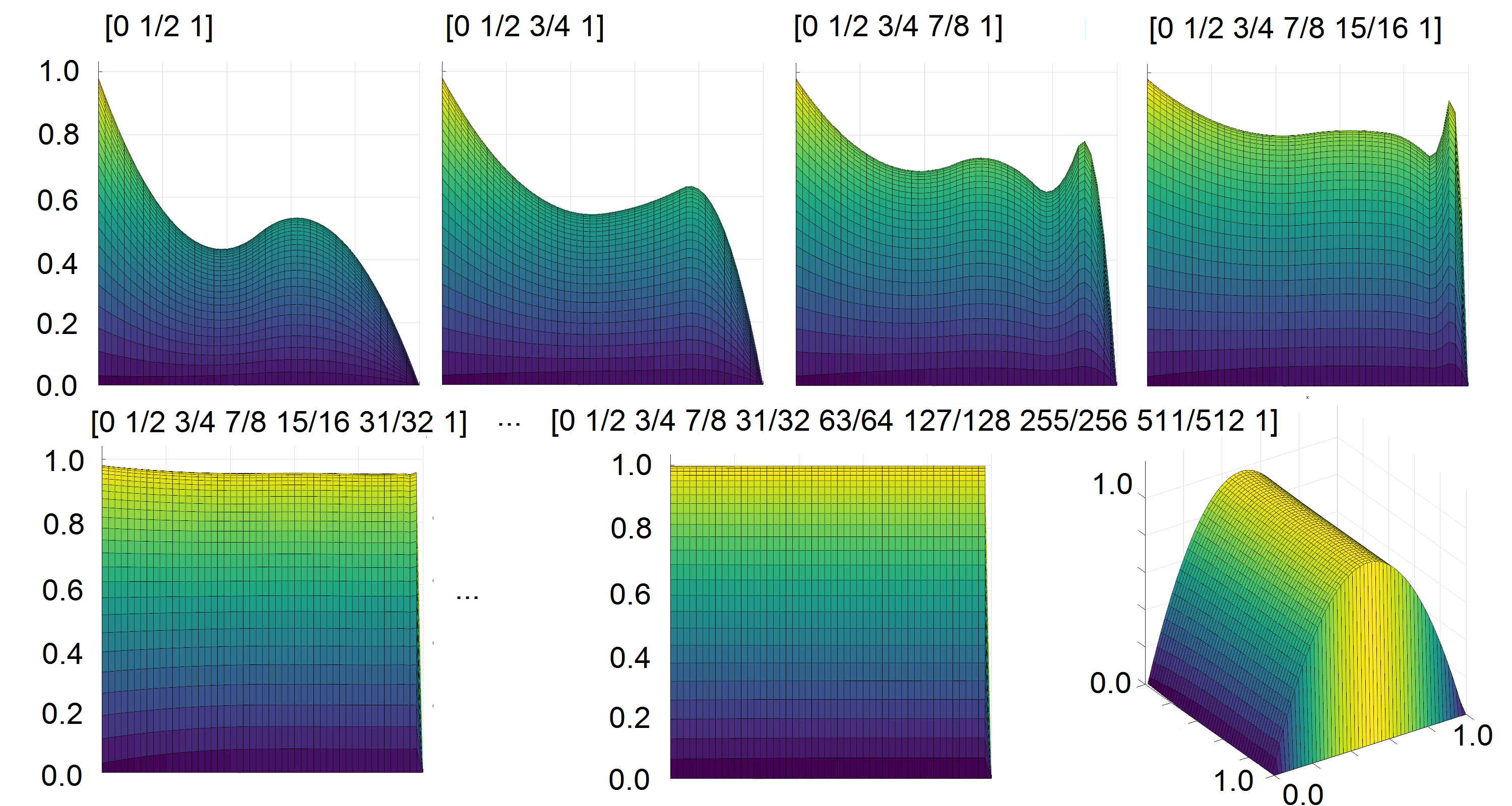}
\caption{Erikkson-Johnsson problem, residual minimization, $\varepsilon=10^{-3}$, using quadratic B-splines $C^1$ for trial and cubic B-splines $C^2$ for testing}
\label{fig:Problem3k}
\end{figure}

To solve the Eriksson-Johnson problem with the finite element method, we need to apply a special stabilization method. We select two methods, the Streamline Upwind Petrov-Galerkin method (SUPG) \cite{HUGHES198797} and the residual minimization method \cite{RM}.
We select $\epsilon=10^{-3}$.
Both methods require adapted mesh. The sequence of solutions from the SUPG method on the adapted mesh is presented in Figure \ref{fig:Problem3l}.
They approximate the solution quite well from the very beginning, but they provide a smooth solution (they do not approximate the stiff gradient well).
The sequence of solutions obtained from the residual minimization method is presented in Figure \ref{fig:Problem3k}. First, they deliver some oscillations, and the SUPG solution is much better on the coarse mesh, but once the mesh recovers the boundary layer, the resulting solution from the Petrov-Galerkin method delivers a good approximation.

\section{Neural network approach for two-dimensional advection-dominated diffusion model problem}

\subsection{Physics Informed Neural Networks for Eriksson-Johnson problem}

We now introduce the PINN formulation of the Eriksson-Johnson problem.
Starting from the strong form: Given $\Omega=(0,1)^2$, $\beta=(1,0)^T$, we seek $\Omega \ni (x,y) \rightarrow u(x,y)$ such that $
\frac{\partial u}{\partial x}-\epsilon \left(\frac{\partial^2 u}{\partial x^2}+
\frac{\partial^2 u}{\partial y^2}\right)=0 \nonumber
$
with Dirichlet boundary conditions
$
u=0 \textrm{ for }x\in(0,1),y\in\{0,1\}$ and $u=sin(\Pi y) \textrm{ for }x=0$,
we introduce the neural network as the solution of the PDE

\begin{equation}
u(x,y)=\inred{NN(x,y)}=A_n \sigma\left(A_{n-1}\sigma\left(...\sigma\left(A_1\begin{bmatrix} x \\ y \end{bmatrix} +B_1\right)...\right)+B_{n-1}\right)+B_n
\end{equation}

We define the following loss functions for the strong form of the PDE, as well as for the Dirichlet boundary conditions:
\begin{eqnarray}
LOSS_{PDE}(x,y) = (\frac{\partial \inred{NN(x,y)}}{\partial x}-\epsilon \frac{\partial^2\inred{NN(x,y)}}{\partial x^2}-\epsilon \frac{\partial^2\inred{NN(x,y)}}{\partial y^2})^2, \\
LOSS_{BC0y}(0,y) = ( \inred{NN(0,y)-\sin(\Pi y)})^2, \\
LOSS_{BC1y}(1,y) = ( \inred{NN(1,y)})^2, \\
LOSS_{BCx0}(x,0) = ( \inred{NN(x,0)})^2,
LOSS_{BCx1}(x,1) = (\inred{NN(x,1)})^2, \\
LOSS(x,y) = LOSS_{PDE}(x,y)+LOSS_{BC0y}(0,y)+LOSS_{BC1y}(1,y)+ \notag \\ +LOSS_{BCx0}(x,0)+LOSS_{BCx1}(x,1)
\end{eqnarray}
The argument for the loss function is the points $(x,y)\in (0,1)^2$, selected during the training process.

\subsection{Variational Physics Informed Neural Networks for Eriksson-Johnson problem}

{We start from the weak formulation: Find $u\in H^1([0,1]^2)$: $b(u,v)=l(v)$
\begin{eqnarray}
b(u,v)=\int_{[0,1]^2} \epsilon\frac{\partial \inblack{u}}{\partial x}\frac{\partial \inbrown{v}}{\partial x}dxdy+
\int_{[0,1]^2} \epsilon\frac{\partial \inblack{u}}{\partial y}\frac{\partial \inbrown{v}}{\partial y}dxdy+
\int_{[0,1]^2}\frac{\partial \inblack{u}}{\partial
x}\inbrown{ v}dxdy
\\ l(v)=0
\end{eqnarray}
For test functions we take cubic B-splines such as $B_{ij;3}(x,y)=B^x_{i,3}(x)B^y_{j,3}(y)$, $i=1,...,N_x$, $j=1,...,N_y$ on $\partial \Omega$.
Now, the neural network represents the PDE solution
\begin{equation}
u(x,y)=\inred{NN(x,y)}=A^n \sigma\left(A^{n-1}\sigma(...\sigma(A^1\begin{bmatrix} x \\ y \end{bmatrix}+b^1)...+b^{n-1}\right)+b^n
\end{equation}
where $A^k_{ij}$ are the trainable neural network weights from layer $k$, and $b^k_i$ are the trainable coefficients of neural network biases.
We define the first loss function by averaging the strong form with test functions
\begin{eqnarray}
LOSS_{strong}= \notag \\ \left(\int_{[0,1]^2} \left( \frac{\partial \inred{NN(x,y)}}{\partial x}-\epsilon \left(\frac{\partial^2 \inred{NN(x,y)}}{\partial x^2}+
\frac{\partial^2 \inred{NN(x,y)}}{\partial y^2}\right)\right) \inbrown{ B^xB^y} dxdy \right)^2,
\end{eqnarray}
where the test functions $B^xB^y$ are cubic B-splines.
We also define the second loss based on the weak formulation:
\begin{eqnarray}
LOSS_{weak}(\gamma) = \left(
b(\inred{NN(x,y)},\gamma \cdot B^xB^y)-l(\inred{NN(x,y)}, B^xB^y)\right)^2.
\end{eqnarray}
We also introduce the loss functions for the Dirichlet boundary conditions:
\begin{eqnarray}
LOSS_{BC}(x,y)=( \inred{NN(0,y)-\sin(\Pi y)})^2+ ( \inred{NN(1,y)})^2+ \notag \\ ( \inred{NN(x,0)})^2+(\inred{NN(x,1)})^2.
\end{eqnarray}
We use the Adam optimization algorithm \cite{kingma2017adam}.}

\subsection{Numerical results}

\subsubsection{PINN and VPINN on uniformly distributed points}
The first numerical experiment concerns the PINN method with points $(x,y)$ selected for the training as a uniform distribution, as well as the VPINN method with cubic B-splines employed as test functions span over the uniformly distributed points.
We employ 4 layers with 20 neurons per layer, hyperbolic tangent as the activation function, and the learning rate is defined as $\eta = 0.00125$.
We use Adam optimizer \cite{kingma2017adam}.
We executed the experiments for different values of $\varepsilon=0.1,0.01,0.001$ with different numbers of points (or intervales), $X=50\times 50,80\times 80$.
PINN and VPINN with different loss functions do not provide correct numerical results when using uniform distribution of points or test functions.

\subsubsection{PINN and VPINN on adapted mesh}

\begin{figure}[h]
\includegraphics[scale=0.7]{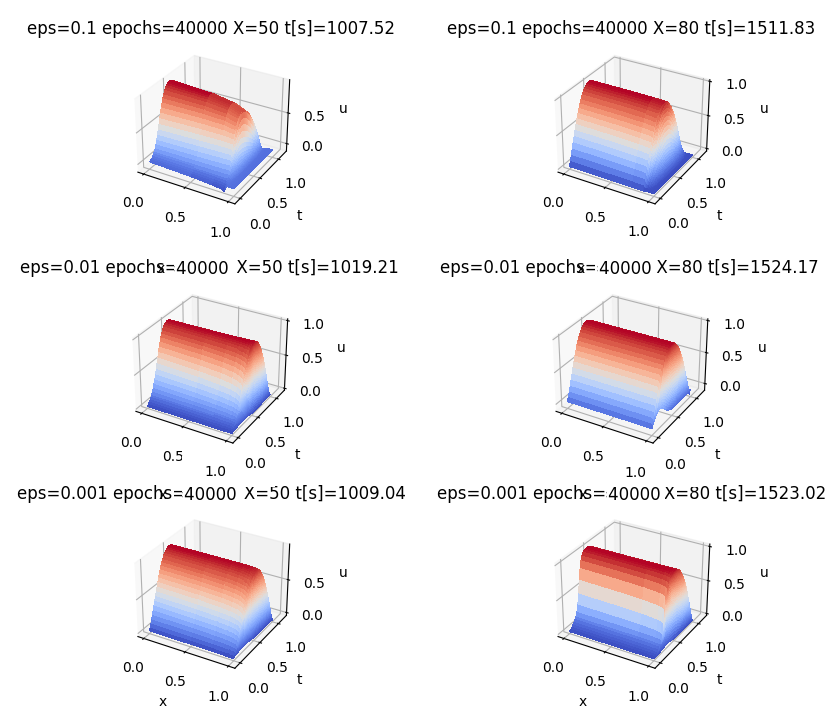}
\caption{Numerical results for PINN method for Eriksson-Johnson problem on adapted mesh}
\label{fig:results_adaptive_EJ_PINN}
\end{figure}

The second numerical experiment concerns the PINN method with points $x$ selected for the training based on the adapted mesh and with cubic B-splines selected for the VPINN span on the adapted mesh.
The mesh has been defined as $x_0=0, x_1=0.5, x_i=x_{i-1}+\frac{x_{i-1}+x_{i-2}}{2}$, up to the point where $1-x_i<\epsilon$ and then we put equally distributed remainning points between $1-\epsilon$ and 1.
We employ 4 layers with 20 neurons per layer, hyperbolic tangent as the activation function, and the learning rate is defined as $\eta = 0.00125$.
The PINN results are summarized in Figure \ref{fig:results_adaptive_EJ_PINN}.
The rows correspond to different values of $\varepsilon=0.1,0.01,0.001$ and columns correspond to different number of points (or intervales), $X=100,1000$.
We can read from these Figures that using adapted mesh, it is possible to solve the Eriksson-Johnson problem using the PINN method with arbitrary $\epsilon=0.1,0.01$, and $0.001$ using 50 or more points.
We can conclude that for 50 adaptively distributed points, the PINN method provides a good solution for $\varepsilon=0.01$. There are 40,000 epochs, and the total training time is 1019 seconds.
We can also conclude that for 50 adaptively distributed points, the PINN method provides a good solution for $\varepsilon=0.001$. There are 40,000 epochs, and the total training time is 1009 seconds.

\begin{figure}[h]
\includegraphics[scale=0.7]{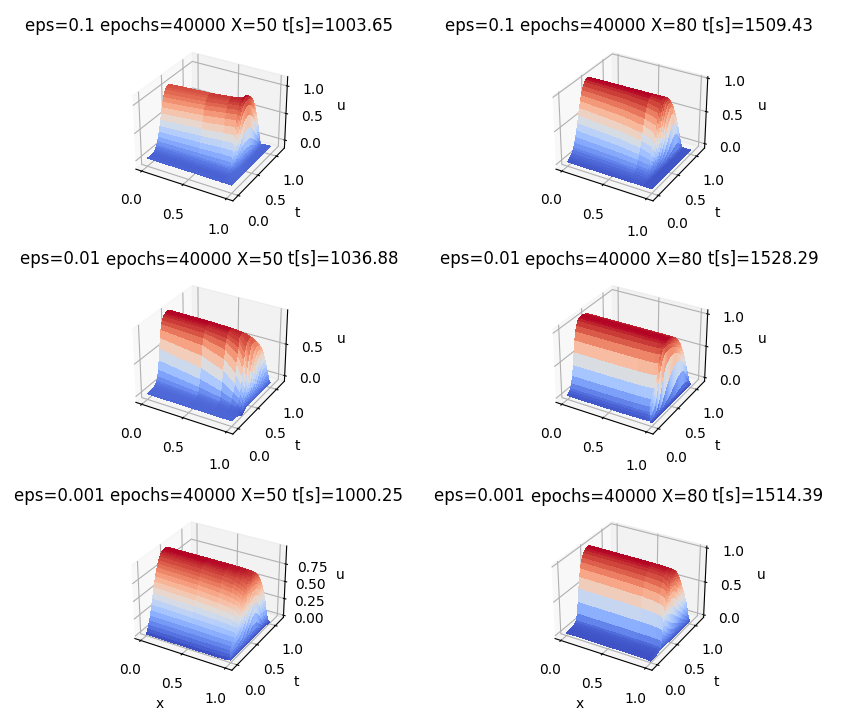}
\caption{Numerical results for VPINN with strong loss function for Eriksson-Johnson problem on adapted mesh}
\label{fig:results_adaptive_EJ_VPINN_strong}
\end{figure}

\begin{figure}[h]
\includegraphics[scale=0.7]{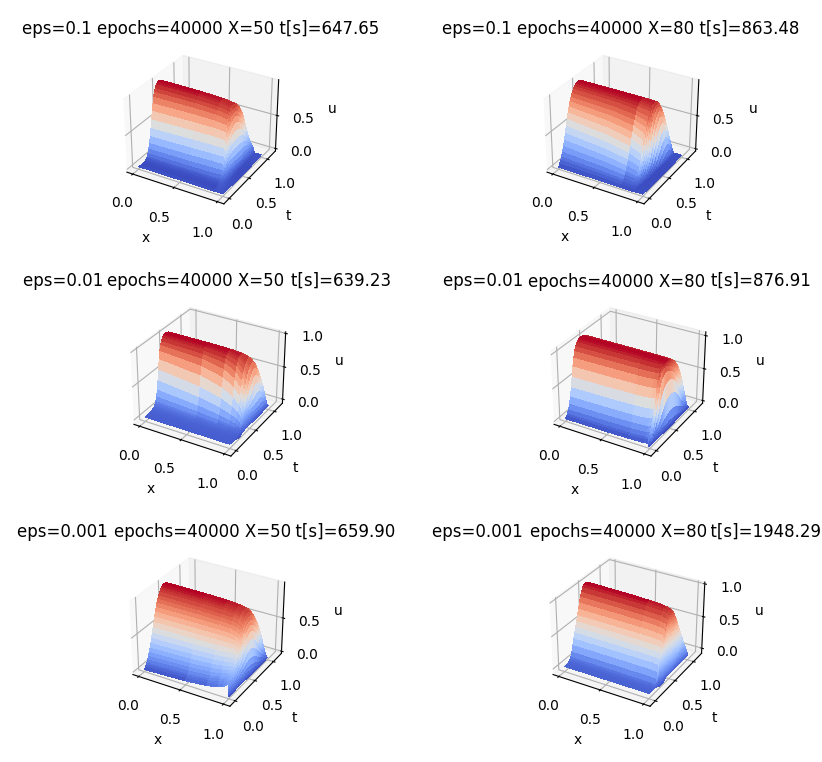}
\caption{Numerical results for VPINN with weak loss function for Eriksson-Johnson problem on adapted mesh}
\label{fig:results_adaptive_EJ_VPINN_weak}
\end{figure}

\begin{figure}[h]
\includegraphics[scale=0.7]{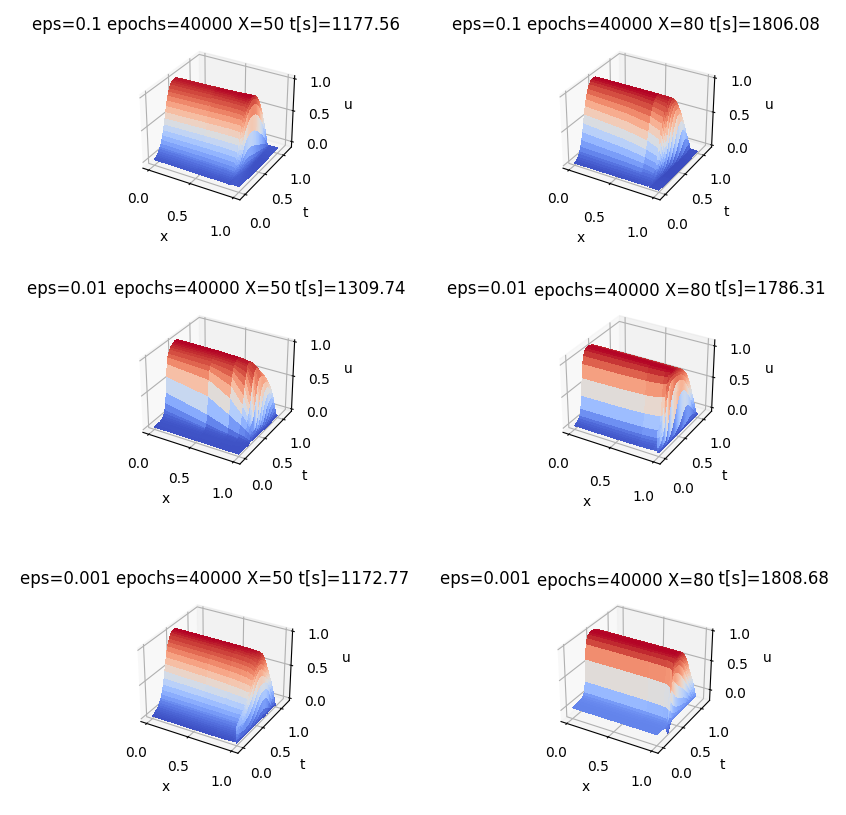}
\caption{Numerical results for VPINN with weak and strong loss functions together for Eriksson-Johnson problem on adapted mesh}
\label{fig:results_adaptive_EJ_VPINN_weak_strong}
\end{figure}

The VPINN results with strong loss are summarized in Figure \ref{fig:results_adaptive_EJ_VPINN_strong}.
The VPINN results with weak loss are summarized in Figure \ref{fig:results_adaptive_EJ_VPINN_weak}.
The VPINN results with weak and strong loss functions together are summarized in Figure \ref{fig:results_adaptive_EJ_VPINN_weak_strong}.

The rows correspond to different values of $\varepsilon=0.1,0.01,0.001$, and columns correspond to a different number of points (or intervales), $X=100,1000$.
We can read from these Figures that if using adapted mesh, it is possible to solve the Eriksson-Johnson problem using the VPINN method with arbitrary loss functions and arbitrary $\epsilon=0.1,0.01$, and $0.001$ using 50 or more points.
Tables \ref{tab:PINN} and \ref{tab:VPINN} illustrate the numerical accuracy with respect to the exact solution
\begin{eqnarray}
u_{exact}(x, y) = \frac{(e^{(r_1 (x-1))} - e^{(r_2 (x-1))})}{ (e^{(-r_1)} - e^{(-r_2)})} sin(\Pi y), \\ r_1 = \frac{(1 + \sqrt{(1 + 4\epsilon^2\Pi^2)})}{ (2\epsilon)}, r_2 = \frac{(1 - \sqrt{(1 + 4\epsilon^2\Pi^2)})}{ (2\epsilon)}.
\end{eqnarray}

\begin{center}
\begin{table}
\begin{tabular}{ | c |c | c|}
\hline
Number of points & $\epsilon$ & MSE$(u,u_{exact})$ \\
\hline
X=50 & $\epsilon=0.1$ & MSE$(u,u_{exact})$=0.05\\
X=80 & $\epsilon=0.1$ & MSE$(u,u_{exact})$=0.08 \\
\hline
X=50 & $\epsilon=0.01$ & MSE$(u,u_{exact})$=0.18\\
X=80 & $\epsilon=0.01$ & MSE$(u,u_{exact})$=0.2 \\
\hline
X=50 & $\epsilon=0.001$ & MSE$(u,u_{exact})$=0.15\\
X=80 & $\epsilon=0.001$ & MSE$(u,u_{exact})$=0.2\\
\hline
\end{tabular}
\caption{Numerical accuracy of PINN solution of the Eriksson-Johson problem}
\label{tab:PINN}
\end{table}
\end{center}

\begin{center}
\begin{table}
\begin{tabular}{ | c |c | c|}
\hline
Number of points & $\epsilon$ & MSE$(u,u_{exact})$ \\
\hline
X=50 & $\epsilon=0.1$ & MSE$(u,u_{exact})$=0.09\\
X=80 & $\epsilon=0.1$ & MSE$(u,u_{exact})$=0.07 \\
\hline
X=50 & $\epsilon=0.01$ & MSE$(u,u_{exact})$=0.06\\
X=80 & $\epsilon=0.01$ & MSE$(u,u_{exact})$=0.12 \\
\hline
X=50 & $\epsilon=0.001$ & MSE$(u,u_{exact})$=0.11\\
X=80 & $\epsilon=0.001$ & MSE$(u,u_{exact})$=0.11\\
\hline
\end{tabular}
\caption{Numerical accuracy of VPINN solution with weak and strong loss functions of the Eriksson-Johson problem}
\label{tab:VPINN}
\end{table}
\end{center}

\begin{figure}[h]
\includegraphics[scale=0.45]{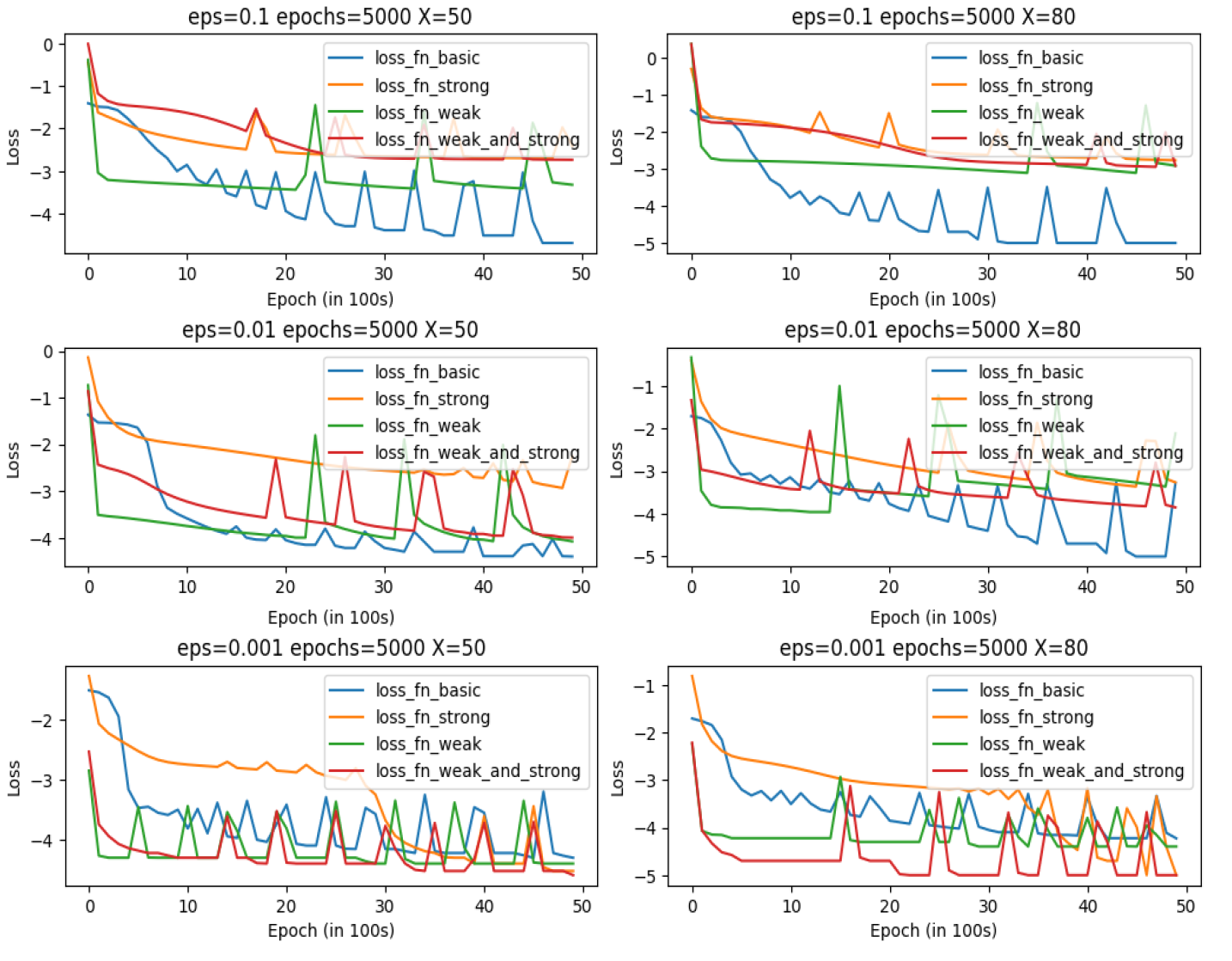}
\caption{Convergence of PINN and VPINN training for Eriksson-Johnson on adapted mesh}
\label{fig:results_uniform_EJ}
\end{figure}

\section{Conclusions}

The neural network can be successfully applied for the solution of challenging advection-dominated diffusion problems.
We focused on Physics Informed Neural Networks (PINN) using the strong residual as the loss function,
as well as on the Variational Physics Informed Neural Networks (VPINN) using strong and weak residuals (multiplied by test functions and possibly integrated by parts).

The training algorithm selects points to evaluate the residual loss for PINN during the training. It also selects test functions to span over the selected mesh. The points and test functions are selected to evaluate the residuals for VPINN during the training process.

When using a uniform distribution of 100 points or a uniform mesh of 100 intervals to span the test functions, it is only possible to train VPINN with strong residual to find the solution of the model one-dimensional problem with $\epsilon=0.01$. Using 1000 uniform points, or 1000 elements mesh, does not allow finding the correct solution for $\epsilon=0.001$, even using 150,000 epochs.
For the two-dimensional Eriksson-Johnson problem, the uniform distribution of points does not work at all.
On the other hand, when using adapted mesh, both PINN and VPINN provide the correct solution for the one-dimensional model problem and two-dimensional Eriksson-Johnson problem.

Thie higher-order finite element methods behave similarly. The adaptive mesh is a must to solve the two-dimensional Eriksson-Johnson problem using stabilized weak formulations.
From the practical point of view, finite element method solvers require the development of stabilized formulations, like Streamline Upwind Petrov-Galerkin, or residual minimization methods, together with the weak imposition of the boundary conditions \cite{CALO2021113214}.

The (Variational) Physics Informed Neural Networks employ strong or weak residual and strong enforcement of the boundary conditions. They work fine with the diffusion coefficient $\epsilon=0.001$. The Adam optimization algorithm finds a good-quality solution after 40,000 epochs. The training time takes around 1000-2000 seconds on a modern laptop. For smaller values of diffusion, e.g., $\epsilon \leq 0.0001$, better loss functions related to the Riesz representative of the residual may be required, which will be a subject of our future work.

Artificial intelligence methods based on PINN/VPINN technology, are attractive alternatives for the solution of challenging engineering problems.

One disadvantage of the method is the need for pre-existing knowledge about the location of the boundary layer to construct an efficient adaptive mesh.
Our future work will include the development of automatic adaptive algorithms for training PINN/VPINN methods.

\section*{Acknowledgments}

Research project partly supported by the program "Excellence Initiative – research university" for the AGH University of Science and Technology.


\end{document}